# Direct linearization method for nonlinear PDE's and the related kernel RBFs

W. Chen

Department of Informatics, University of Oslo, P.O.Box 1080, Blindern, 0316 Oslo, Norway

Email: wenc@ifi.uio.no

## Abstract

The standard methodology handling nonlinear PDE's involves the two steps: numerical discretization to get a set of nonlinear algebraic equations, and then the application of the Newton iterative linearization technique or its variants to solve the nonlinear algebraic systems. Here we present an alternative strategy called direct linearization method (DLM). The DLM discretization algebraic equations of nonlinear PDE's is simply linear rather than nonlinear. The basic idea behind the DLM is that we see a nonlinear term as a new independent systematic variable and transfer a nonlinear PDE into a linear PDE with more than one independent variable. It is stressed that the DLM strategy can be applied combining any existing numerical discretization techniques. The resulting linear discretization equations can be either over-posed or well-posed. In particular, we also discuss how to create proper radial basis functions in conjunction with the DLM.

**Key words**: Nonlinear PDE's**,** Newton-type methods, direct linearization method, kernel radial basis function.

# 1. Introduction

Although great endeavor has been devoted to nonlinear computation and analysis, it seems very difficult to attack nonlinear problems directly. The linearization procedures such as the Newton method and its variants [1] are now commonly used to transform a nonlinear system to a linear system in a point-wise approximate way so that the standard numerical linear algebra approach can be employed for computation and analysis. It is noted that the strategy of linearization often leads to a very huge amount of computing effort and encounters great difficulty in nonlinear stability analysis. For instance, a nice initial guess solution is key to guarantee the convergence and reliability of the solution, which, however, is often a daunting task. Consequently, a question is aroused: is iteration truly unavoidable? The present study is a new quest among the others [2,3] to attack this problem in a different manner.

The basic idea behind this study is that we see a nonlinear term as a new independent variable and transfer a nonlinear PDE of an independent variable into a linear PDE with more than one independent variable. Then we can apply any standard numerical discretization technique to analogize this linear PDE. To get the well-posed or over-posed discretization formulations, we need to use staggered nodes a few times more of what the standard method requires. It turns out that the size of the resulting system matrix grows linearly proportional to the number of nonlinear terms. This methodology is called the direct linearization method (DLM). It is stressed that the DLM strategy can be applied with any existing numerical discretization techniques and essentially eliminates iteration linearization computation. In conjunction with the DLM for nonlinear PDE's this paper also discusses three kernel RBF strategies.

## 2. Direct linearization method

In order to illustrate our idea clearly, let us consider the quadratic nonlinear equation without loss of generailty

$$p(u)q(u) + \Re(u) = f(x), \qquad (1)$$

where $p(u)$, $q(u)$ and $\Re(u)$ are linear differential operators, $f(x)$ is inhomogeneous term. The mathematical description of the problem is complimented with Dirichlet and Neumann boundary conditions

$$u = \bar{u}, \text{ on } \Gamma_1 \qquad (2a)$$

$$q = \partial u/\partial n = \bar{q}, \text{ on } \Gamma_2 \qquad (2b)$$

where $n$ is the outward normal to the boundary, $\Gamma = \Gamma_1 + \Gamma_2$, and the upper bars indicate known boundary values.

The Newton-type methods of point-wise iterative linearization are standard technique to solve nonlinear analog equations of this PDE system. Instead, if we regard the quadratic nonlinear term as a new independent variable, i.e.,

$$v = p(u)r(u). \qquad (3)$$

In this way, we rewrite nonlinear equation (1) as a linear equation with two independent variables $u(x)$ and $v(x)$

$$v + \Re(u) = f(x). \qquad (4)$$

Accordingly, we have the Dirichlet boundary conditions for the new system independent variable *v*. It, however, is not easy to get Neumann boundary conditions of *v*. We have two approaches to solve this problem. The first is to use the particular solution method [4] or dual reciprocity method [5], which splits the original solution into linear homogeneous solution and nonlinear inhomogeneous particular solution. For the nonlinear particular solution part, we do not need to satisfy the boundary conditions at all. Therefore, we circumvent the Neumann boundary condition issue of the DLM. The second approach is that we simply replace Neumann boundary conditions of the DLM independent variable *v* with the governing equation. Since the original Neumann boundary conditions have been satisfied by variable *u*, this approximation may not any introduce significant errors.

The above steps clarify the basic idea behind the direct linearization method. Now we can use any numerical discretization technique to analogize linear equation (4). For instance, variable *u* and *v* are approximately represented respectively by a finite series approximation

$$u(x) = \sum_{k=1}^{N} \alpha_k \varphi_k(x), \qquad (5)$$

$$v(x) = \sum_{k=1}^{N} \beta_k \phi_k(x), \qquad (6)$$

where $\alpha$ and $\beta$ are unknown expansion coefficients, $\phi$ and $\varphi$ are the basis functions. N is the total number of boundary and inside-domain points. Now we need to discuss how to choose basis functions. Is there any mutual constrains in choosing $\phi$ and $\varphi$ once one of them is determined? The similar issue also appears in handling a set of linear partial differential equations with more than one independent variable. We think that $\phi$ and $\varphi$ should satisfy somehow mutual constrains. However, for this moment, this is still an open question and probably problem-dependent. In section 3, we will give a brief discussion on this issue related to the radial basis function (RBF) approach.

Since we have two independent variables $u$ and $v$ and one equation (4), we need to evaluate 2N unknown coefficients. To get a well-posed or over-posed discretization system of equations, we need to discretize Eqs. (4), (2a,b), and boundary conditions of $v$ at least at 2N points. Namely, we require 2N staggered field points at boundary and inside domain. Then we have 2N discretization equations

$$A_1\alpha + A_2\beta = b_1 \tag{7a}$$

$$C_1\alpha + C_2\beta = b_2 \tag{7b}$$

where $A_1$, $A_2$, $C_1$ and $C_2$ are interpolation matrices of order $N$. Note that formulations (7a,b) correspond to two sets of N staggered points across the problem domain and are a set of linear simultaneous equations. Therefore, no iterative linearization technique is required to solve the well-posed discretization equations (7a,b). If we discretize equation (4), (2a,b), and boundary conditions of $v$ at more than 2N points, we get an over-posed linear system of equations. Then a linear least square method should be used to solve the DLM discretization equations. For more complicated coupled nonlinear PDE systems, it is very straightforward to apply the DLM.

## 3. Constructing kernel RBFs for the DLM

In most practically significant cases, the nonlinear terms in the PDE system is a coupling of several linear operators. The construction of nonlinear algorithm should consider this nonlinear feature. Following this idea, this section discusses the construction of radial basis functions relating to the DLM.

Chen and Tanaka [6,7] points out an underlying relationship between the RBF approximation and Green identity. A kernel-RBF strategy is proposed accordingly. By Green's second theorem, we have solution of Eqs. (1)-(2a,b)

$$u(z) = \int_\Omega w^*_\Re(z,x)[-p(u)q(u) + f(x)]d\Omega + \int_\Gamma \left[ u \frac{\partial w^*_\Re(z,x)}{\partial n} - w^*_\Re(z,x) \frac{\partial u}{\partial n} \right] d\Gamma, \quad (8)$$

where $w^*_\Re$ denotes the fundamental solution of operator $\Re\{\}$, $x$ indicates source point. The above formula (8) suggests us that the exponential augmented kernel (EAK) RBF can be created by

$$\omega(r,x) = [f(x) + p(u)q(u)] r^{2m} w^*_\Re(r), \quad (9)$$

where $m$ is an integral number and $r^{2m}$ augmented term enhances the smoothness and ensures sufficient degree of differential continuity since the fundamental solution has a singularity at origin.

It is worth pointing out here that we can use the nonsingular general solution instead of singular fundamental solution in various kernel RBFs presented in this paper. For the brevity, we only mention the fundamental solution relating to the kernel RBF from now on.

In the RBF approximation, the influence coefficients are only-point dependent. Therefore, one factor decisive to the efficiency is to choose approximate influence functions, which is termed as the RBF normally. It is practically impossible for us to find accurate nonlinear influence functions for complex problems. Simply removing the nonlinear term in the RBF (9), we have an operational EAK RBF

$$\omega(r,x) = f(x) r^{2m} w_{\Re}^{*}(r). \tag{10}$$

The above RBF can be employed to analogize differential systems (1) and (2a,b) with the standard numerical procedure. On the other hand, in terms of the DLM, we need to have two RBFs to respectively approximate two independent variables $u$ and $v$ in Eq. (4).

$$u = \sum_{k=1}^{N} A_k \psi_u(r_k) \tag{11}$$

$$v = \sum_{k=1}^{N} B_k \psi_v(r_k) \tag{12}$$

Here a natural choice for $\psi_u$ is the kernel RBF for linear operator $\Re\{\}$

$$\psi_u = r^{2m} w_{\Re}^{*}(r), \tag{13}$$

while $\psi_v$ should reflect nonlinear feature of operator $p(u)q(u)$, i.e.,

$$\psi_v(r) = p(u)q(u) r^{2n} w_{\Re}^{*}(r). \tag{14}$$

The key issue here is how to simplify the RBF (14). One solution is to use the fundamental solution of operators $p(u)$ and $q(u)$. Therefore, we have

$$\psi_v(r) = r^{2n} w_{\Re}^*(r) w_p^*(r) w_q^*(r). \tag{15}$$

It is quite clear here that we use different operator-dependent kernel RBFs to approximate independent variable $u$ and $v$, since they have different physical meanings. From physical field viewpoint, we reason that the RBF is in fact to evaluate influence coefficient of a source point in terms of influence function of a physical problem. The RBF (15) reflects to some extent physical background of nonlinear terms. In other words, the present scheme is an intrinsically nonlinear numerical discretization technique.

The broad definition of the kernel RBF involves the use of nonsingular high-order fundamental solution and general solution of operators and shape parameter. The second approach creating kernel RBFs for $u$ and $v$ approximation is

$$\psi_u = w_{\Re}^{*m}(r), \tag{16}$$

$$\psi_v(r) = w_{\Re}^{*n}(r) w_p^{*n}(r) w_q^{*n}(r), \tag{17}$$

where $m$ and $n$ are the order of nonsingular high-order fundamental solutions. This strategy is called high-order fundamental solution kernel (HSK) RBF.

By analogy with using the shape parameter in the MQ RBF, we have the third kernel RBF creating strategy:

$$\psi_u = w_{\Re}^*\left(\sqrt{r^2 + c^2}\right), \tag{18}$$

$$\psi_v(r) = w^*_{\Re}\left(\sqrt{r^2+c^2}\right) w^*_p\left(\sqrt{r^2+c^2}\right) w^*_q\left(\sqrt{r^2+c^2}\right), \tag{19}$$

where $c$ is the shape parameter. Substituting $\sqrt{r^2+c^2}$ into the RBFs (16) and (17) instead of $r$ is also an alternate approach. We simply name the present approach of the shape parameter kernel (SPK) RBF.

All complete fundamental solutions consist of essential and complementary elementary functions [8]. The standard singular fundamental solutions used in the BEM involve only the essential part. The complementary terms of the complete fundamental solutions are often understood the nonsingular general solution in terms of the boundary knot method [6,7]. The shape parameter $c$ in the SPK and MQ RBFs can be interpreted as the scaling parameter in the simplified form of the complete fundamental solutions and leads to infinite smoothness. More precisely, the SPK methodology constructs the RBF by substituting $\sqrt{r^2+c^2}$ into the essential terms of the complete fundamental solutions [6,7] to compromise the complementary terms (general solutions). For instance, the reciprocal MQ is physically formed by this methodology using the fundamental solutions of more than 3 dimensional Laplace operators. The MQ is based on the fundamental solutions of 1D Laplacian. The SPK thin plate splines are also presented in [6,7]. Computational benefits using the SPK depend on the tricky choose of the shape parameter, which is in agreement with the skillful implementation of general fundamental solutions. The other possible approach embedding the complementary terms into the essential terms of the complete fundamental solution is still open.

In many cases, we may have no fundamental solutions of operator $p(u)$ and $q(u)$. The following scheme may be a cheap alternative for creating RBF for $v$:

$$\psi_v(r) = r^{2s} p(\psi_u) q(\psi_u), \tag{20}$$

where *s* is an integral number. $\psi_v$ can be understood approximate nonlinear influence function, relative to conventional linear influence function. The essential philosophy behind the present RBF for the DLM is that we should think more the RBF from a physical field point of view than form a mathematical interpolation.

Even if we use the standard single RBF for solving nonlinear differential system (1) and (2a,b) without relating to the DLM, the chosen RBF should be dependent on all differential operators. Based on the fundamental solution of linear operator $\Re$, the EAK RBF is

$$\omega(r, x) = \left[ -p(w_\Re^*(r)) q(w_\Re^*(r)) + f(x) \right] r^{2m} w_\Re^*(r). \tag{21}$$

The HSK RBF is

$$\omega(r, x) = \left[ -p(w_\Re^{*m}(r)) q(w_\Re^{*m}(r)) + f(x) \right] w_\Re^{*m}(r). \tag{22}$$

The SPK RBF is

$$\omega(r, x) = \left[ -p(w_\Re^*(\sqrt{r^2 + c^2})) q(w_\Re^*(\sqrt{r^2 + c^2})) + f(x) \right] w_\Re^*(\sqrt{r^2 + c^2}). \tag{23}$$

## 4. Some remarks

The basic ideas behind the preceding three kernel RBF strategies for the DLM are also applicable to the polynomial approximation and general nonlinear data processing. All in all, we should understand nonlinear solver from system physical essence [8]. This study is still in a very early stage. The practical numerical experiments of this strategy will be provided subsequently.